\documentclass{amsart}

\newtheorem{theorem}{Theorem}[section]
\newtheorem{lemma}[theorem]{Lemma}
\newtheorem{cor}[theorem]{Corollary}
\theoremstyle{definition}
\newtheorem{definition}[theorem]{Definition}

\theoremstyle{remark}

\newcommand{\sN}{\mathcal{N}} 
\newcommand{\sC}{\mathcal{C}} 

\DeclareMathOperator{\End}{End}
\DeclareMathOperator{\im}{im}        

\newcommand{\sneq}{\mbox{$\underline{\triangleleft\mspace{-1.8mu}\triangleleft} \medspace$}}
\newcommand{\id}{\triangleleft}
\newcommand{\ideq}{\unlhd}
\newcommand{\str}{\overrightarrow}
\usepackage{amsfonts}
\usepackage{amscd}  
\numberwithin{equation}{section}

\begin{document}

\title{Engel subalgebras of $n$-Lie algebras} 
\author{Donald W. Barnes}
\address{1 Little Wonga Rd, Cremorne NSW 2090 Australia}
\email{donwb@iprimus.com.au}
\thanks{This work was done while the author was an Honorary Associate of the
School of Mathematics and Statistics, University of Sydney.}
\subjclass[2000]{Primary 17B05, 17B30}
\keywords{$n$-Lie algebras, soluble, nilpotent, Engel subalgebras}

\begin{abstract} Engel subalgebras of finite-dimensional $n$-Lie algebras are shown to have
similar properties to those of Lie algebras.  Using these, it is shown that an $n$-Lie algebra, all of whose maximal
subalgebras are ideals, is nilpotent.  A primitive $2$-soluble $n$-Lie algebra is shown to split over its minimal ideal   
and that all the complements to its minimal ideal are conjugate.   A subalgebra is shown to be a Cartan subalgebra if  
and only if it is minimal Engel, provided that the field has sufficiently many elements.  Cartan subalgebras are
shown to have a property analogous to intravariance.  
\end{abstract}
\maketitle
\section{Introduction}
In this section, I set out the basic definitions and notations used.  A fuller account of the basic theory of $n$-Lie algebras
is given in Kasymov \cite{Kas}.  In Section $2$, I set out the basic properties of soluble and nilpotent algebras and of
their representations.  In Section $3$, I establish the properties of Engel subalgebras and use them to prove some
analogues of known theorems on Lie algebras.  In Section $4$, I show that, provided that the field has at least $\dim(L) +
1$ elements, the Cartan subalgebras of the $n$-Lie algebra $L$ are precisely its minimal Engel subalgebras.  I show that 
Cartan subalgebras have a property which could be regarded as the $n$-Lie analogue of the group theory concept of
intravariance.  (A subgroup $U$ of a group $G$ is called {\em intravariant} in $G$ if every automorphism of $G$ maps
$U$ onto a conjugate subgroup.)

\begin{definition} An $n$-Lie algebra over a field $F$ is a vector space $L$ over $F$ together with a
linear map $m: \Lambda^nL: \to L$, usually written $m(x_1 \wedge \dots \wedge x_n) = [x_1, \dots, x_n]$, satisfying  
the generalised Jacobi identity
\begin{equation}\label{eq1}
[[x_1, \dots, x_n], y_2, \dots, y_n] = \sum_{i=1}^n [x_1, \dots, [x_i, y_2, \dots, y_n], \dots,
x_n].\end{equation}
\end{definition}

For $n = 2$, this is the definition of a Lie algebra.  Note that if $n > 2$ and we take a fixed element $a
\in L$, then the map $m': \Lambda^{n-1} L \to L$ given by $m'(x_1 \wedge \dots \wedge x_{n-1}) =
m(x_1 \wedge \dots \wedge x_{n-1} \wedge a)$ makes $L$ into an $(n-1)$-Lie algebra.  I shall only consider 
finite dimensional algebras.

The left multiplication map $D(a_1, \dots, a_{n-1}): L \to L$ given by  
$$D(a_1, \dots, a_{n-1})(x) = [a_1,  \dots, a_{n-1}, x]$$ 
is clearly a derivation and will be called an inner derivation, as will  any linear combination of 
such derivations.  (Kasymov uses right multiplication.)  To shorten  notation, where it can be 
done without ambiguity, an array $a_1, \dots, a_r$ will be written $\str{a}$, without 
subscripts, whatever the length of the array.  Thus we write
$[\str{x}]$ for $[x_1, \dots, x_n]$ and $D(\str{a})$  for  $D(a_1, \dots, a_{n-1})$.   For a subspace $A$ of $L$ and a string
$\str{a}$, I write $\str{a} \in A$ if all $a_i \in A$ and denote by $D(A)$ the space $\langle D(\str{a}) \mid  \str{a} \in A
\rangle$ spanned by the $D(\str{a})$ with $\str{a} \in A$.   If $A \le L$ then $D(A)$ is a subalgebra of the Lie algebra of
derivations of $L$.

I  introduce the some notations for operations on  subspaces.  I denote by $[A_1, \dots, A_n]$  the
subspace spanned by all the
$[\str{a}]$ with
$a_i \in A_i$.  I put 
$$[r_1A_1, r_2A_2, \dots, r_pA_p] = [A_1, \dots, A_1, A_2, \dots, A_2,\dots, A_p],$$
where there are $r_i$ factors $A_i$, and $\sum_ir_i = n$.  If $r_i = 1$, I omit it from the notation.  

A subspace $K \subseteq L$ is called a subalgebra (written $K \le  L$,) if $[K, \dots K] \subseteq
K$ and is called an ideal (written $K \ideq L$,) if $[K, L, \dots, L] \subseteq K$.  If $K \ideq L$, then 
we can form the quotient algebra $L/K$.   If $A, B \ideq L$,  we put $A \cdot B = [A,B, (n-2)L]$. 

\begin{definition}  An algebra $L$ is called \textit{abelian} if $[nL] = 0$.  An ideal $A \ideq L$ is 
said to be \textit{abelian embedded} in $L$ if  $[2A, (n-2)L] = 0$.\end{definition}

\begin{definition}  An $L$-module is a vector space $V$ together with a map $\rho: 
\Lambda^{n-1}L \to \End(V)$ such that setting $[\str{x}, v] = \rho(\str{x})(v)$ and $[2A,
(n-2)L] = 0$ makes $L \oplus V$ into an $n$-Lie algebra.  The map $\rho$ is called the
representation afforded by the module $V$.
\end{definition}

Equivalently, the map $\rho: \Lambda^{n-1}L \to \End(V)$ is a representation if 
\begin{equation}\label{eq2}
\rho(\str{x}) \rho(\str{y})(v) = \sum_{i=1}^{n-1} \rho(y_1, \dots, [\str{x}, y_i], \dots
y_{n-1})(v)
\quad + \rho(\str{y}) \rho(\str{x}(v),\end{equation}
or equivalently, 
\begin{equation}\label{eq3}
[\rho(\str{x}), \rho(\str{y})] = \sum_{i=1}^{n-1} \rho(y_1, \dots, [\str{x}, y_i], \dots,
y_{n-1}),\end{equation}
and
\begin{equation} \label{eq-lrep2}
\rho([x_1, \dots, x_n], y_2, \dots, y_{n-1}) = \sum_{i=1}^n (-1)^{n-i} \rho(x_1, \dots,
\widehat{x_i}, \dots,  x_n) \rho(x_i, y_2, \dots, y_{n-1}).\end{equation}

Condition (\ref{eq2}) ensures the generalised Jacobi identity for those products where an element of $V$ occurs in the
inner of the two products,     but does not deal with the cases where  the inner multiplication is of $n$ elements of
$L$, with an element of $V$ as a factor in the outer product.

Note that the subspace $\rho(L)$ of $\End(V)$ spanned by all the $\rho(\str{x})$ for $\str{x} 
\in L$ is a Lie  algebra of linear transformation of $V$.  A subspace $W \subseteq V$ is an
$L$-submodule if and only if it is a $\rho(L)$-submodule.

The $n$-Lie algebra $L^*$ which we have defined on the vector space $L \oplus V$ is called the split extension of
$V$ by $L$.  $V$ is an abelian embedded ideal of $L^*$ and $L$ is a subalgebra of $L^*$.  The set $\ker(\rho)
= \{k \in L \mid \rho(k, \str{x}) = 0$ for all $\str{x} \in L\}$ is called the kernel of
the representation $\rho$ or of the module $V$.  It is an ideal of $L$.  In particular, we may regard $L$ as an
$L$-module affording the representation $D$ of inner derivations.  This is called the regular or adjoint 
representation.  Its kernel is the centre $Z(L) = \{z \in L \mid [\str{x}, z] = 0 \text{ for all }\str{x} \in L\}$.  The
inner derivations of $L$ form a Lie algebra.   Any ideal $A \ideq L$ may be regarded as a submodule of the adjoint
module.  The kernel of the representation of $L$ afforded by $A$ is the centraliser 
$$\sC_L(A) = \{c \in L \mid [c, A, (n-2)L] = 0\}.$$
The centraliser of an ideal, being the kernel of a representation, is an ideal.  If the ideal $A$ is abelian embedded 
in $L$, then $A \subseteq  \sC_L(A)$.  Having $A$ merely abelian is not sufficient for this.

\section{Solubility and nilpotency}
There are several  definitions of solubility and nilpotency in the literature.  Rather than use phrases such as   
``soluble in the sense of Filippov'', I use more descriptive terminology.

\begin{definition} An algebra $L$ is called \textit{soluble} if $L^{(r)} = 0$ for some $r$, where $L^{(0)} = L$   
and  $L^{(s+1)} = [nL^{(s)}]$. \end{definition}

\begin{definition} An ideal $K \ideq L$ is said to be \textit{$k$-solubly embedded in $L$}, $2 \le   k  \le n$, if
$K^{(r,k)} = 0$ for some $r$, where $K^{(0,k)} = K$ and $K^{(s+1,k)} = [kK^{(s,k)},
(n-k)L]$. \end{definition}

\begin{definition} An algebra $L$ is called \textit{nilpotent} if $L^r = 0$ for some $r$, where $L^1 = L$ and
$L^{s+1} = [L^s, (n-1)L]$. \end{definition}

\begin{definition} An ideal $K \ideq L$ is said to be \textit{$k$-nilpotently embedded in $L$} if $K^{r,k} = 0$ for 
some $r$, where $K^{1,k} = K$ and $K^{s+1,k} = [K^{s,k}, (k-1)K, (n-k)L]$.
\end{definition}
We have the following results of Kasymov \cite{Kas}.
\begin{lemma} Let $A\supseteq B$ be ideals of $L$.  Suppose $B$ is $k$-solubly embedded in $L$ and that $A/B$
is $k$-solubly embedded in $L/B$.  Then $A$ is $k$-solubly embedded in $L$. \end{lemma}
\begin{cor} Suppose $A,B$ are $k$-solubly embedded ideals of $L$.  Then $A+B$ is a $k$-solubly embedded ideal.
\end{cor}

\begin{definition} The largest $k$-solubly embedded ideal of $L$ is called the \textit{$k$-radical} of $L$ and 
denoted by $S_k(L)$.  (It is the sum of all the $k$-solubly embedded ideals of $L$. \end{definition}

I note some trivial properties.  If $L \ne 0$ is soluble, then $L$ has a non-zero abelian ideal.  If $L \ne 0$ is 2-soluble, 
then $L$ has a non-zero ideal which is abelian embedded in $L$.  If $L \ne 0$ is nilpotent, then $Z(L) \ne 0$.  If $A$
is an abelian embedded ideal of $L$ and $B \subseteq A$ is an ideal of $L$, then $B$ is abelian embedded in $L$. 
Thus a $2$-soluble algebra has a minimal ideal which is abelian embedded.

\begin{lemma}\label{lem-compl}  Let $A$ be a minimal abelian embedded ideal of $L$ and suppose $U < L$ and
that $U+A = L$.  Then $U \cap A = 0$ and $U$ is a maximal subalgebra of $L$.  \end{lemma}

\begin{proof}  Put $B = U \cap A$.  Since $A$ is abelian embedded, $[B, A, (n-2)L] = 0$.  But $L = A + U$, so it
follows that $[B, (n-1)L] = [B, (n-1)U] \subseteq A \cap U$.  Thus $B \id L$ and from the minimality of $A$, it follows
that $B = 0$.  That $U$ is maximal follows. \end{proof}

\begin{definition}  Let $U \le L$.  The \textit{normaliser} $\sN_L(U)$ of $U$ in $L$ is the set 
$$\{x \in L \mid [x, u_1, \dots, u_{n-1}] \in U \text{ for all } u_1, \dots, u_{n-1} \in U\}.$$ \end{definition}

Note that $\sN_L(U)$ is a subspace of $L$.  It need not be a subalgebra.  If $x \in \sN_L(U)$, then the subspace $\langle
x, U\rangle$ spanned by $x$ and $U$ is a subalgebra and $U$ is an ideal of $\langle x, U\rangle$.  If $U \ne L$ and $L$ is
nilpotent, then 
$U \ne \sN_L(U)$.  It follows that, if $L$ is nilpotent, then every maximal subalgebra of $L$ is an ideal.  

\begin{lemma}\label{lem-sums} Let $L$ be a nilpotent, $n$-Lie algebra. Then every inner derivation of $L$ is nilpotent. 
\end{lemma}
\begin{proof}  For every string $\str{x} \in L$, we have $D(\str{x}) L^k \subseteq L^{k+1}$.  Thus for all $d \in D(L)$, we
have $d L^k \subseteq L^{k+1}$. \end{proof}

This has a converse, the analogue of Engel's Theorem, due to Kasymov
\cite[Theorem 3]{Kas}.

\begin{theorem}\label{th-Eng} Suppose $D(\str{x})$  is nilpotent for all $\str{x}$ in  the finite dimensional 
$n$-Lie algebra $L$.  Then $L$ is nilpotent. \end{theorem}
This is an immediate consequence of Kasymov \cite[Theorem 4]{Kas}.

\begin{theorem}\label{th-linEng} Suppose $\rho$ is a faithful representation of the $n$-Lie
algebra $L$ and that all the  $\rho(\str{x})$ are nilpotent.  Let $A$ be the associative algebra
generated by the $\rho(\str{x})$.  Then $A$  and $L$ are nilpotent. \end{theorem}

For a subset $U \subseteq L$ and representation $\rho$, I put $\rho_r(U) = \langle \rho( \str{x})
\mid x_1, \dots, x_r \in U \rangle$ and write simply $\rho(U)$ for $\rho_{n-1}(U)$.   The next two lemmas are Kasymov
\cite[Proposition 3]{Kas}.

\begin{lemma}\label{lem-nilrep}  Let $\rho$ be a representation of the $n$-Lie
algebra $L$.  Then 
$$\rho_1(L^k) \subseteq \rho_1(L)^k.$$  If $L$ is nilpotent, then $\rho_1(L)$ is nilpotent.
\end{lemma}

\begin{lemma}  Let $S$ be a $2$-solubly embedded ideal of $L$ and let $\rho$ be a 
representation of $L$.  Then $\rho_1(S)$ is soluble. \end{lemma}

\begin{definition} A subalgebra $U \le L$ is called \textit{subnormal} in $L$, written $U \sneq L$, if there exists a 
chain of subalgebras $U_0 = U < U_1 < \dots < U_r = L$ with each $U_i \id U_{i+1}$. \end{definition}

\begin{definition} Let $N \sneq L$.  Let $A$ be an irreducible $N$-module.  An $L$-submodule
$W$ of the $L$-module $V$ is called an $A$-component of $V$ if every $N$-composition factor 
of $W$ is isomorphic to $A$ while no $N$-composition factor of $V/W$ is isomorphic to $A$.
\end{definition}

\begin{theorem}\label{th-sort} Let $L$ be an $n$-Lie algebra and let $V$ be an $L$-module. Suppose that $N \sneq L$   
is nilpotent.  Then for each irreducible $N$-module $A_i$, there exists an $A_i$-component $A_i(V)$ and $V = \oplus_i
A_i(V)$.  \end{theorem}

\begin{proof} Let $\rho$ be the representation afforded by $V$.  Note that $\rho(N) \sneq \rho(L)$ and is
nilpotent.  The result follows by Barnes \cite{sort}.
\end{proof}
 
\begin{cor}  Let $N$ be a nilpotent ideal of $L$ and let $V$ be an irreducible $L$-module.  Then
all $N$-composition factors of $V$ are isomorphic.\end{cor}

If $N \ideq L$ and $\rho$ is the representation of $L$ afforded by the module $V$, then $\rho_k(N) \ideq \rho(L)$.    If
$N$ is $k$-nilpotently embedded in $L$, then $\rho_k(N)$ is nilpotent.  The $\rho_k(N)$ components of $V$
are $\rho(L)$-submodules and  so are also $L$-submodules.  This gives a possibly finer direct decomposition of $V$. 
These direct decompositions are natural. 

\section{Engel subalgebras}
Engel subalgebras of Lie algebras, so named because of their obvious connection to Engel's Theorem, have the 
useful property that any subalgebra containing one is self-normalising.  I show that the analogously defined
subalgebras of an $n$-Lie algebra have the analogous property.  An inner derivation of a Lie algebra is given by a single
element.  For an $n$-Lie algebra, an inner derivation need not be given by a single string, but is a sum of derivations
$D(\str{x})$ given by strings $\str{x}$.  To cope with this, we have to put an extra condition into the definition of Engel
subalgebras to make them work as they do for Lie algebras.  Let $d$ be any derivation of $L$.  We put
$$E_L(d) = \{x \in L \mid d^rx = 0 \text{ for some }r\}.$$ 

\begin{lemma}  Let $d$ be any derivation of $L$.  Then $E_L(d) \le L$. \end{lemma}
\begin{proof}    Let $e_1, \dots, e_n \in E_L(d)$.  For each $i$, we have $r_i$ such
that $d^{r_i}e_i = 0$.  Put $r = \sum_i r_i$.   As $d$ is a derivation, $d^r[\str{e}]$ is a sum 
of terms of the form $[d^{s_1}e_1, \dots, d^{s_n}e_n]$ with $\sum_i s_i = r$.  For some $j$,
we must have $s_j \ge r_j$, so  the term is $0$.  Thus $[e_1, \dots, e_n] \in E_L(d)$.  \end{proof}

\begin{definition} Let $L$ be an $n$-Lie algebra and let $d$ be an inner derivation of $L$.  
The subalgebra $E = E_L(d)$ is called an Engel subalgebra if $d \in D(E)$.\end{definition}
Note that, for any string $\str{a} \in L$, we have $D(\str{a})a_i =  0$, so $\str{a} \in E_L(D(\str{a}))$ and 
$E_L(D(\str{a}))$ is an Engel subalgebra.

\begin{lemma}\label{lem-selfN} Let $E = E_L(d)$ be an Engel subalgebra of $L$.  If $E \le U \le  L$, then
$\sN_L(U) = U$.
\end{lemma}

\begin{proof}  Let $N = \sN_L(U)$.   For any $\str{e} \in E$, we have $D(\str{e})N \subseteq U$, and so $dN \subseteq
U$.   Let $d_N = d| N \to N$ and let $r = \dim(N)$.  Then $N = \im(d_N^r) \oplus \ker(d_N^r) \subseteq U + E = U$. 
\end{proof}

\begin{theorem}\label{th-WB}  Suppose every maximal subalgebra of the $n$-Lie algebra $L$ is an ideal.  Then
$L$ is nilpotent. \end{theorem}

\begin{proof}  Let $D = D(\str{x})$  and let $E = E_L(D)$ be the corresponding 
Engel subalgebra.  If $E \ne L$, there exists a maximal subalgebra $M$ of $L$ containing $E$. 
But then $\sN_L(M) = M$ contrary to $M$ being an ideal.  Therefore every $D(\str{x})$ is
nilpotent and $L$ is nilpotent by Theorem \ref{th-Eng}. \end{proof}

\begin{definition} The \textit{Frattini subalgebra} $\Phi(L)$ is the intersection of the maximal subalgebras of $L$.
\end{definition}

Theorem \ref{th-WB} can be stated in terms of the Frattini subalgebra:  $L$ is nilpotent if and only if $\Phi(L)
\supseteq L \cdot L$.  If $L$ is nilpotent, then in fact $\Phi(L) = L \cdot L$.

Bai has proved \cite[Theorem 2.4]{Bai} that any ideal of $L$ contained in $\Phi(L)$ is nilpotent.  The following two
theorems strengthen and generalise that result.

\begin{theorem}\label{th-GB} Let $U \sneq L$, $K \ideq U$ with $K \subseteq \Phi(L)$.  Suppose that $U/K$ is
nilpotent.  Then $U$ is nilpotent. \end{theorem}

\begin{proof}  We have a chain of subalgebras $U_0 = U \id U_1 \id \dots \id U_r = L$.  Let $\str{a}\in U$ and 
put $D = D(\str{a})$.  Then $DU_i \subseteq U_{i-1}$ since $U_{i-1} \id U_i$.  Hence $D^rL
\subseteq U$.  But $U/K$ is nilpotent, so $D^sU \subseteq K$ for some $s$.   Thus, if $t = \dim(L)$, we have 
$D^tL \subseteq K$.  But $L = \im(D^t) \oplus \ker(D^t)$, so $L = K + E_L(D)$.  But $K
\subseteq \Phi(L)$, so this implies that $E_L(D) = L$.  Thus every $D(\str{a})$ for $\str{a} \in U$ is
nilpotent and $U$ is nilpotent by Theorem \ref{th-Eng}. \end{proof}

\begin{theorem} \label{th-Gasch2nil} Let $A \subseteq B$ be ideals of the $n$-Lie algebra $L$. 
Suppose that $B/A$ is $2$-nilpotently embedded in $L/A$ and that $A \subseteq \Phi(L)$.  Then
$B$ is $2$-nilpotently embedded in $L$. \end{theorem}

\begin{proof} Let $\str{x} \in L$ have $x_1 \in B$.  Since $B \ideq L$, $D(\str{x})L \subseteq B$.  Since $B/A$ is
$2$-nilpotently embedded in $L/A$, $D^r(\str{x})B \subseteq A$ for some $r$.   Thus $D(\str{x})^rL
\subseteq A$ for some $r$ and we have $E_L(\str{x}) + A = L$.  As $A \subseteq \Phi(L)$, we have
$E_L(\str{x}) = L$ for all such $\str{x}$.   Therefore $B$ is $2$-nilpotently embedded in $L$ by Theorem \ref{th-linEng}.
\end{proof}

A soluble Lie algebra with a self-centralising minimal ideal is called primitive.  In Barnes \cite[Theorem 4]{Cosol} and
Barnes and Newell \cite[Theorem 1.1]{BN}, it was proved that a primitive Lie algebra splits over its minimal ideal and
that all complements are conjugate.  The analogous result holds for $2$-soluble $n$-Lie algebras.  We first need the
appropriate notion of conjugacy.  Let $d: L \to L$ be a derivation of $L$ such that $d^2=0$ and $[2(dL),(n-2)L] = 0$. 
Then
$\alpha = 1 + d$ is an automorphism of $L$ since
\begin{equation*} \begin{split}
[(1+d)x_1, \dots, (1+d)x_n] = &[x_1, \dots, x_n] + \sum_{i=1}^n [x_1, \dots, dx_i, \dots, x_n]\\
 & \quad +\text{ terms of higher degree in $d$}.\end{split}\end{equation*}
The terms of degree greater than 1 in $d$ are, by assumption, all zero.  Thus $[\alpha x_1, \dots, \alpha x_n] =
\alpha [x_1, \dots, x_n]$.   In particular, if $A$ is an abelian embedded ideal of $L$, then any inner derivation
$D(a,\str{x})$ with $a \in A$ satisfies the conditions.  The corresponding automorphism $\alpha$
will  be called an $A$-inner automorphism as will any product of such.  Note that any two such
commute since $A$ is  abelian embedded. If $U, V \le L$ and $\alpha U = V$ for some $A$-inner
automorphism $\alpha$, we say that
$U$ and $V$ are $A$-conjugate. The following lemma deals with an awkward point in the proof. 

\begin{lemma}\label{lem-nil} Let $A$ be a minimal abelian embedded ideal of $L$ and let $B/A$ be an abelian
embedded ideal of $L/A$.  For $b \in B$ and $x_1, \dots, x_{n-2} \in L$, put 
$$\rho(b, \str{x}) = D(b, \str{x})|A \to A.$$  
Suppose every $\rho(b, \str{x})$ is nilpotent.  Then  $[A, B, (n-2)L] = 0$.
\end{lemma}
 
\begin{proof}  By Theorem \ref{th-linEng}, the associative algebra generated by the $\rho(b, \str{x})$ acts nilpotently   
on $A$.  Thus $[A, B, (n-2)L] \ne A$.  But $[A, B, (n-2)L] \ideq L$.  As $A$ is a minimal ideal of $L$, it follows that $[A, B,
(n-2)L] = 0$. \end{proof}

\begin{theorem}\label{th-prim} Let $L$ be a 2-soluble $n$-Lie algebra and let $A$ be a minimal abelian embedded
ideal.  Suppose $\sC_L(A) = A$.  Then $L$ splits over $A$ and all complements to $A$ in $L$ are $A$-conjugate.
\end{theorem}

\begin{proof}  If $L = A$, then the result holds, so we may suppose $L \ne A$.  Take $B \ideq L$ such that $B/A$ is a
minimal abelian embedded ideal of $L/A$.  Let $\rho$ be the representation of $L$ afforded by $A$.  By asumption,
we have $\ker(\rho) = A$.  If every $\rho(b, \str{x})$ with $b \in B$ is nilpotent, then $[A,B, (n-2)L]
= 0$ by Lemma \ref{lem-nil}, that is, $B \subseteq \sC_L(A)$ contrary to hypothesis.  Thus there exists $b, 
\str{x}$ such that   $d = D(b, \str{x})$ is not nilpotent.  We have $dL \subseteq B$ and
$d^2L \subseteq A$ since $B/A$ is abelian embedded in $L/A$.  For $t = \dim(L)$, we have $L = \ker(d^t)
\oplus \im(d^t)$ and $\im(d^t) \subseteq A$.  Thus $A + E_L(b, \str{x}) = L$ but $ E_L(b, \str{x}) \ne L$ and 
$U = E_L(b, \str{x})$ complements $A$ in $L$.

Now let $U = E_L( \str{x})$ be a complement to $A$ in $L$ with $x_1 \in B$, and let $V$ be another
complement.  Let $y_i = x_i + a_i$ be the element of $V$ in the coset $x_i + A$.  We have from above that $d =
D(\str{x})$ acts invertibly on $A = \im(d^t)$.  Hence there exists $a'_i \in A$ such that $da'_i = a_i$.  Put
$\alpha_i = 1 + (-1)^{n-i}D(x_1, \dots, \hat{x}_i, \dots, x_{n-1}, a'_i)$.  For $j \ne i$, $D(x_1, \dots, \hat{x}_i,
\dots, x_{n-1}, a'_i)(x_j) = [x_1, \dots, \hat{x}_i, \dots, x_{n-1}, a'_i, x_j] = 0$ because of the repeated $x_j$.
Thus $\alpha_i(x_i) = y_i$ and $\alpha_i(x_j) = x_j$ for $j \ne i$.  Put $\alpha = \Pi_i \alpha_i$.  Then $a(x_i) =
y_i$ for all $i$.  The isomorphism $U \to L/A \to V$ maps $x_i$ to $y_i$.  $D(\str{x})$ acts nilpotently on $U$, so
$D(\str{y})$ acts nilpotently on $V$ and $E_L(\str{y}) \supseteq V$.  But $E_L(\str{y}) = \alpha E_L(\str{x})$, so
has dimension $\dim(U) = \dim(V)$.  Thus $E_L(\str{y}) = V$ and $V = \alpha U$.  
\end{proof}

In a primitive Lie algebra in which the quotient by the minimal ideal is abelian, any two complements 
have zero intersection.  For the $n$-Lie algebra above with $n > 2$, we have $x_j \in U \cap \alpha_iU$ for $j \ne i$
even if $L/A$ is abelian.  However, the intersection of all the  complements is zero.  If $u \in U$, $u \ne 0$, then $u
\not\in
\sC_L(A)$, so there exist $\str{y}
\in L$ and $a \in A$ such that $[\str{y}, u, a] = a' \ne 0$.  Put $\theta = 1 + D(\str{y}, a)$.  Then $\theta u = u -a' \in \theta
U$.  If also $u \in
\theta U$, we would have $a' \in A \cap \theta U = 0$.

Theorem \ref{th-prim} opens up the possibility of a theory of formations and projectors for
$2$-soluble $n$-Lie algebras, analogous to the theory developed by Gasch\"utz for finite soluble groups in \cite{Gasc}
and set out in Doerk and Hawkes \cite{DH}.  The Lie algebra version is given in Barnes and Gastineau-Hills  \cite{BGH}.
\section{Cartan subalgebras}
As for Lie algebras, we have 
\begin{definition} A Cartan subalgebra of an $n$-Lie algebra is a nilpotent subalgebra $S$ such
that $\sN_L(S) = S$. \end{definition}

These have been studied in Kasymov \cite{Kas}.
\begin{lemma}\label{lem-minE} Let $L$ be an\ \ $n$-Lie algebra of dimension $k$ over a field F of
at least $k+1$ elements.  Let $U \le L$ and put $D(U) = \langle D(\str{x}) \mid \str{x} \in U
\rangle$.  Suppose that $E = E_L(d_0)$ is minimal in the set $\{E_L(d) \mid d \in D(U)\}$ and
that $E \supseteq U$.  Then $E_L(d) \supseteq E$ for all $d \in D(U)$. \end{lemma}

\begin{proof}  $L$ may be regarded as a $U$-module under the action of the $d \in D(U)$.  
Since $E \le L$ and $E \supseteq U$, $E$ is a $U$-submodule.  Take any $d_1 \in D(U)$ and let
$\theta(t, \lambda), \phi(t, \lambda)$ and $\psi(t, \lambda)$ be the characteristic polynomials
of  $d_0 + \lambda d_1$ on $L, E$ and $L/E$ respectively.  Then
\begin{equation*} \begin{split}
\theta(t, \lambda) &= \phi((t, \lambda) \psi(t, \lambda)\\
\phi(t, \lambda) &=t^r + \alpha_1(\lambda) t^{r-1} + \dots + \alpha_r(\lambda)\\
\psi(t, \lambda) &= t^{k-r} + \beta_1(\lambda) t^{k-r-1} + \dots + \beta_{k-r}(\lambda)\\
\end{split}\end{equation*}
where $r = \dim(E)$ and $\alpha_i(\lambda), \beta_i(\lambda)$ are polynomials in $\lambda$
of degree at most $i$.  We prove $\alpha_i = 0$ for all $i$.

Since $0$ is not an eigenvalue of $d_0$ on $L/E$, 
$\beta_{k-r}(0) \ne 0$.  Thus $\beta_{k-r}(\lambda)$ is not the zero polynomial.  Since
$\beta_{k-r}(\lambda)$ has at most $k-r$ roots in $F$, there exists $r+1$ distinct elements
$\lambda_1, \dots, \lambda_{r+1} \in F$ such that $\beta_{k-r}(\lambda_j) \ne 0$.  But
$\beta_{k-r}(\lambda_j) \ne 0$ implies that $E_L(d_0 + \lambda_j d_1) \subseteq E$.  Therefore
$E_L(d_0 + \lambda_j d_1) = E$ by the minimality of $E$.  But this implies $\phi(t, \lambda_j)
= t^r$, that is, $\alpha_i(\lambda_j) = 0$ for all $j$.  Since $\alpha_i$ has degree at most $i <
r+1$, $\alpha_i(\lambda)$ is the zero polynomial.  Thus $E_L(d_0 + \lambda d_1) \supseteq E$
for all $d_1 \in D(U)$ and all $\lambda \in F$.  Given $d \in D(U)$, put $d_1 = d_0- d$.  Then
$E_L(d) = E_L(d_0 + d_1) \supseteq E$.
\end{proof}

\begin{theorem}\label{th-minEng} Let $L$ be an $n$-Lie algebra of dimension $k$ over the 
field $F$ of at least $k+1$ elements.  Let $S \le L$.  Then $S$ is a Cartan subalgebra of $L$ if 
and only if $S$ is minimal in the set of Engel subalgebras of $L$. \end{theorem}

\begin{proof}  Suppose $S = E_L(d)$ is minimal Engel in $L$.  Then $\sN_L(S) = S$ by Lemma
\ref{lem-selfN}.  For any $\str{s} \in S$, $E_L(\str{s}) \supseteq S$ by Lemma 
\ref{lem-minE}.  Thus all  inner derivations of $S$ are nilpotent and $S$ is nilpotent by 
Theorem \ref{th-Eng} and $S$ is a Cartan subalgebra.

Suppose conversely, that $S$ is a Cartan subalgebra.  Then for any $d \in D(S)$, we have
$E_L(d) \supseteq S$ by Lemma \ref{lem-sums}, so  $S$ does not properly contain any Engel subalgebra of $L$.    We
have to show that $S$ is an Engel subalgebra of $L$, that is, that there exists $d \in D(S)$ such that 
$E_L(d) = S$.  Take $d_0 \in D(S)$ such that $E = E_L(d_0)$ is minimal in the set $\{E_L(d)
\mid d \in D(S)\}$.  By Lemma \ref{lem-minE}, $E_L(d) \supseteq E$ for all $d \in D(S)$.  Thus
$S$ is represented on the $S$-module $E/S$ by nilpotent linear transformations.  If $E/S \ne
0$, then by Theorem \ref{lem-nilrep}, there exists a $1$-dimensional submodule $T/S$ of
$E/S$.  The representation of $S$ on $T/S$ is $0$, that is $D(\str{s})T \subseteq S$ for all
$\str{s} \in S$.  Thus $T \subseteq \sN_L(S) = S$.
\end{proof}

\begin{cor} Let $L$ be an $n$-Lie algebra of dimension $k$ over the field $F$ of at least $k+1$
elements.  Then $L$ has a Cartan subalgebra. \end{cor}
This strengthens Kasymov \cite[Theorem 5]{Kas}.

A subgroup $U$ of a group $G$ is called {\em intravariant} if, for every automorphism $\alpha$ of $G$, $\alpha(U)$ is
conjugate to $U$ in $G$, or equivalently, if every automorphism of $G$ is the product of an inner automorphism and  
an automorphism which stabilises $U$.  For Lie algebras, the use of derivations seems more appropriate.  A subalgebra
$U$ of the Lie algebra $L$ is called {\em intravariant} if every derivation of $L$ is the sum of an inner derivation and a
derivation which stabilises $U$.  In Barnes \cite[Lemma 1.2]{Frat}, it was shown that a subalgebra $U$ of the Lie algebra
$K$ is intravariant in $K$ if and only if, for every Lie algebra $L$ containing $K$ as an ideal, $L = K + \sN_L(U)$.  For
$n$-Lie algebras, this argument fails.  We cannot extend an $n$-Lie algebra by an arbitrary derivation, and an inner
derivation need not be given by a single string.  However, we do have an analogue of Barnes \cite[Theorem 2.1]{Frat}. 

\begin{theorem}\label{th-intrav}  Let $K \ideq L$ and let $S$ be a Cartan subalgebra of $K$. 
Then $L = K + \sN_L(S)$. \end{theorem}

\begin{proof}   Put $N = \sN_L(S)$ and consider $L/N$ as an $S$-module.  Suppose $x+N \in
L/N$ gives the $1$-dimensional zero representation of $S$.  Then $D(\str{s})x \in N$ for all
$\str{s} \in S$.  But $D(\str{s})x \in   K$ since $\str{s} \in K$ and $K \id L$.  Thus $\str{s}x
\in K \cap N = \sN_K(S) = S$.  Thus $x \in N$.  Threfore $L/N$ has no $1$-dimensional
submodule affording the zero representation.  By Theorem \ref{th-sort}, the $S$-module $L/N$
has no composition factor on which $S$ acts trivially.  But $S$ acts trivially on every composition
factor above $K$ of the $S$-module $L$.  Therefore there are no $S$-composition factors above
both $K$ and $N$.  Therefore $K+N=L$.
\end{proof}

\bibliographystyle{amsplain}

\end{document}